# Demuškin Groups with Operators

by Kay Wingberg at Heidelberg

In this paper we consider pro-$p$ Poincaré groups of dimension 2, also called "Demuškin groups". These groups are defined as follows: a pro-$p$-group $G$ is called a Demuškin group if its cohomology has the following properties:

$$\dim_{\mathbb{F}_p} H^1(G,\mathbb{Z}/p\mathbb{Z}) < \infty,$$

$$\dim_{\mathbb{F}_p} H^2(G,\mathbb{Z}/p\mathbb{Z}) = 1, \text{ and the cup-product}$$

$$H^1(G,\mathbb{Z}/p\mathbb{Z}) \times H^1(G,\mathbb{Z}/p\mathbb{Z}) \xrightarrow{\cup} H^2(G,\mathbb{Z}/p\mathbb{Z}) \text{ is non-degenerate.}$$

In the following we exclude the exceptional case that $G \cong \mathbb{Z}/2\mathbb{Z}$. The dualizing module $I$ of $G$ is isomorphic to $\mathbb{Q}_p/\mathbb{Z}_p$ as an abelian group and we have a canonical action of $G$ on $I$.

Demuškin groups occur as Galois groups of the maximal $p$-extension of $p$-adic number fields (if these fields contain the group of $p$-th roots of unity) and as the $p$-completion of the fundamental group of a compact oriented Riemann surface. In the first case the action of $G$ on $I$ is non-trivial whereas in the second case $G$ acts trivially on $I$. We will only consider Demuškin groups acting non-trivially on its dualizing module and we are interested in free pro-$p$-quotients of these groups. Possible ranks of such free quotients were first calculated in [7], [6] and [2].

In many cases of interest a finite group $\Delta$ of order prime to $p$ acts on a Demuškin group. As an example consider the local field $k = \mathbb{Q}_p(\zeta_p)$, where $p$ is an odd prime number. Then $G(k|\mathbb{Q}_p) \cong \mathbb{Z}/(p-1)\mathbb{Z}$ acts on the Demuškin group $D = G(k(p)|k)$, where $k(p)$ is the maximal $p$-extension of $k$. Of particular interest is the case where $\Delta$ is generated by an involution, e.g. $G(k|\mathbb{Q}_p(\zeta_p + \zeta_p^{-1})) \cong \mathbb{Z}/2\mathbb{Z}$ acts on $D$; see [8] where Demuškin groups with involution were considered.

In the following, we are mainly interested in free pro-$p$-quotients of a Demuškin group $G$ which are invariant under a given action of $\Delta$ on $G$. We show the existence of $\Delta$-invariant free quotients $F$ of $G$ such that $F^{ab}$ has a prescribed action of $\Delta$.

If $p$ is odd and $\varDelta \cong \mathbb{Z}/2\mathbb{Z}$ acts on a $p$-Demuškin group $G$ of rank $n+2$, then there exists a $\varDelta$-invariant free quotient $F$ of $G$ such that $\operatorname{rank}_{\mathbb{Z}_p}(F^{ab})^+ = 1$ and $\operatorname{rank}_{\mathbb{Z}_p}(F^{ab})^- = n/2$. This situation occurs for number fields as the following example shows:

Let $p$ be an odd regular prime number and consider the CM-field $k = \mathbb{Q}(\zeta_p)$ with maximal totally real subfield $k^+ = \mathbb{Q}(\zeta_p + \zeta_p^{-1})$. Then the Galois group $G(k|k+) \cong \mathbb{Z}/2\mathbb{Z}$ acts on the Galois group $G(k_{S_p}(p)|k)$ of the maximal $p$-extension $k_{S_p}(p)$ of $k$ which is unramified outside $p$. Let $k_\mathfrak{p}$ be the completion of $k$ with respect to the unique prime $\mathfrak{p}$ of $k$ above $p$ and let $k_\mathfrak{p}(p)$ its maximal $p$-extension. Then we have a surjection

$$G(k_\mathfrak{p}(p)|k_\mathfrak{p}) \qquad G(k_{S_p}(p)|k)$$

of the Demuškin group $G(k_\mathfrak{p}(p)|k_\mathfrak{p})$ of rank $p+1$ onto the free pro-$p$-group $G(k_{S_p}(p)|k)$ of rank $(p+1)/2$ and

$$\operatorname{rank}_{\mathbb{Z}_p}(G(k_{S_p}(p)|k)^{ab})^+ = 1, \quad \operatorname{rank}_{\mathbb{Z}_p}(G(k_{S_p}(p)|k)^{ab})^- = (p-1)/2\,.$$

It would be of interest under which conditions there exist large free quotients of $G(k_{S_p}(p)|k)$ for an arbitrary CM-field $k$. If we assume that no prime $\mathfrak{p}$ above $p$ splits in the extension $k|k^+$, such a quotient should be defined by free quotients of the local groups $G(k_\mathfrak{p}(p)|k_\mathfrak{p})$, $\mathfrak{p}|p$, with an acting of $G(k_\mathfrak{p}|k_p^+) \cong G(k|k^+)$ as above.

# 1 Operators

Let $p$ be a prime number and let

$$1 \qquad G \qquad \mathcal{G} \overset{s}{\qquad} \varDelta \qquad 1\,,$$

be a split exact sequence of profinite groups where $G$ is a pro-$p$-group and $\varDelta$ is a finite group of order prime to $p$. Thus $\mathcal{G}$ is the semi-direct product of $\varDelta$ by $G$ and $G$ is a pro-$p$-$\varDelta$ operator group where the action of $\varDelta$ on $G$ is defined via the splitting $s$. Conversely, given a pro-$p$-$\varDelta$ operator group $G$, we get a semi-direct product $\mathcal{G} = G \rtimes \varDelta$ where the action of $\varDelta$ on $G$ is the given one.

Let $\mathcal{G}(p)$ be the maximal pro-$p$-quotient of $\mathcal{G}$ and let $G_\varDelta$ be the maximal quotient of $G$ with trivial $\varDelta$-action. Observe that $G_\varDelta$ is well-defined.

**Proposition 1.1** *With the notation and assumptions as above there is a canonical isomorphism*

$$G_\varDelta \xrightarrow{\sim} \mathcal{G}(p)\,.$$

*Furthermore, if $H^2(G,\mathbb{Z}/p\mathbb{Z})^\varDelta = 0$, then $G_\varDelta$ is a free pro-$p$-group.*

**Proof:** Consider the exact commutative diagram

$$\begin{array}{ccccccccc} & & & & & s & & & \\ 1 & & N\cap G & & N & & \Delta & & 1 \\ & & & & & s & & & \\ 1 & & G & & \mathcal{G} & & \Delta & & 1 \\ & & & & & & & & \\ & & \tilde{G} & & \mathcal{G}(p)\,, & & & & \end{array}$$

where $N$ is the kernel of the canonical surjection $\mathcal{G} \twoheadrightarrow \mathcal{G}(p)$ and $\tilde{G}$ denotes the quotient $G/N\cap G$. Since $\Delta$ acts on $N\cap G$ via $s$, we obtain an induced action on $\tilde{G}$. This action is trivial because

$$g^{s(\sigma)-1} = [s(\sigma), g] \in N\cap G \quad \text{for } g\in G \text{ and } \sigma \in \Delta\,,$$

and so we get a surjection

$$\varphi : G_\Delta \qquad \tilde{G}\,.$$

Consider the exact commutative diagram

$$\begin{array}{cccccccc} 0 & H^1(\tilde{G}) & & H^1(G_\Delta) & H^1(\ker\varphi)^{G_\Delta} & H^2(\tilde{G}) & & H^2(G_\Delta) \\ & & & \mathit{inf} & & & & \mathit{inf} \\ & \mathit{res} & & H^1(G)^\Delta & & \mathit{res} & & H^2(G)^\Delta \\ & & & \mathit{res} & & & & \mathit{res} \\ & H^1(\mathcal{G}(p)) & \mathit{inf}_1 & H^1(\mathcal{G}) & & H^2(\mathcal{G}(p)) & \mathit{inf}_2 & H^2(\mathcal{G}) \end{array}$$

where $H^i(-) = H^i(-,\mathbb{Z}/p\mathbb{Z})$, and the bijectivity of $\mathit{inf}_1$ and the injectivity of $\mathit{inf}_2$ follows from $\mathrm{Hom}(N,\mathbb{Z}/p\mathbb{Z}) = 0$. We see that $H^1(\ker\varphi)^{G_\Delta} = 0$, and so $\ker\varphi = 1$, i.e. $G_\Delta \cong \tilde{G} \cong \mathcal{G}(p)$.

Furthermore, it follows that

$$H^2(G_\Delta) \xrightarrow{\mathit{inf}} H^2(G)^\Delta$$

is injective. Therefore, if $H^2(G)^\Delta = 0$, then $H^2(G_\Delta) = 0$, and so $G_\Delta$ is a free pro-$p$-group. □

In the following lemma let $G$ be a finitely generated pro-$p$-$\Delta$ operator group, i.e. there is a homomorphism of profinite groups $\Delta \longrightarrow \mathrm{Aut}(G)$, cf. [4] 5.1. Let

$$1 \longrightarrow R \longrightarrow F \longrightarrow G \longrightarrow 1$$

be a minimal presentation of the group $G$ by a free pro-$p$-group $F$ of rank $\dim_{\mathbb{F}_p} H^1(G,\mathbb{Z}/p\mathbb{Z})$ and a normal subgroup $R$.

**Lemma 1.2** *With the notation as above there exists a continuous action of $\Delta$ on $F$ extending the action on $G$, i.e. the surjection $F \twoheadrightarrow G$ is $\Delta$-invariant.*

**Proof:** We may assume that the action of $\Delta$ on $G$ is faithful, i.e. $\Delta$ injects into $\mathrm{Aut}(G)$. Let $\sigma \in \Delta$. Since $F$ is free, there exists a homomorphism $\tilde{\sigma} : F \to F$ such that the diagram

$$\begin{array}{ccccc} & & & & F \\ & & \tilde{\sigma} & & \\ F & & G & \sigma & G \end{array}$$

commutes. The map $\tilde{\sigma}$ is necessarily an automorphism of $F$ since it induces a bijection on $F/F^2 \cong G/G^2$.

Let $\tilde{\Delta}$ be the subgroup of $\mathrm{Aut}(F)$ generated by all lifting $\tilde{\sigma}$ of the elements $\sigma \in \Delta$. Recall that the kernels of the homomorphisms $\mathrm{Aut}(G) \to \mathrm{Aut}(G/G^2)$ and $\mathrm{Aut}(F) \twoheadrightarrow \mathrm{Aut}(F/F^2)$ are pro-$p$-groups, cf. [4] 5.5, and that $\Delta$ has an order prime to $p$. Therefore, we can consider $\Delta$ as a subgroup of $\mathrm{Aut}(G/G^2) \cong \mathrm{Aut}(F/F^2)$, and there exists a subgroup $\Delta'$ of $\tilde{\Delta}$ such the diagram

$$\begin{array}{ccc} \Delta' \subseteq \tilde{\Delta} & & \mathrm{Aut}(F) \\ & & \\ \Delta & & \mathrm{Aut}(F/F^2) \end{array}$$

commutes. Since every lifting $\tilde{\sigma}$ respects $R$, the same is true for the elements of $\Delta'$, i.e. $\Delta'$ consists of liftings of elements of $\Delta$. This proves the lemma. □

Now we assume that

$\Delta$ is a finite group of order prime to $p$ and
$G$ is a $p$-Demuškin group of rank $n+2, n \geq 0$, with dualizing module $I$ and an action by $\Delta$.

Let $\mathcal{G}$ be the semi-direct product of $\Delta$ by $G$, i.e. the sequence

$$1 \longrightarrow G \longrightarrow \mathcal{G} \longrightarrow \Delta \longrightarrow 1$$

is exact.

The dualizing module $I$ of $G$ is defined as

$$I = \varinjlim_{m} \varinjlim_{U} H^2(U, \mathbb{Z}/p^m\mathbb{Z})^*,$$

where $U$ runs through the open normal subgroups of $G$ and the second limit is taken over the maps *cor* $^*$ dual to the corestriction; the first limit is taken with respect to the multiplication by $p$.

Let

$$\chi : G \longrightarrow \mathrm{Aut}(I) \cong \mathbb{Z}_p^\times$$

be the character given by the action of $G$ on $I$. We denote the canonical quotient $G/\ker(\chi)$ by $\Gamma$, i.e.

$$\chi_0 : \Gamma \qquad \mathrm{Aut}(I)\,.$$

In the following we assume that

$$G \quad \text{acts non-trivially on} \quad I$$

(thus $\Gamma \cong \mathbb{Z}_p$), and we define the (finite) invariant $q$ of $G$ by

$$q = \#(I^G).$$

Then we have a $\Delta$-invariant isomorphism

$$H^2(G, \mathbb{Z}/q\mathbb{Z}) \cong \mathrm{Hom}(I^G, \mathbb{Z}/q\mathbb{Z}) \quad (\cong \mathbb{Z}/q\mathbb{Z} \text{ as an abelian group})$$

and a $\Delta$-invariant non-degenerate pairing

$$H^1(G, \mathbb{Z}/q\mathbb{Z}) \times H^1(G, \mathbb{Z}/q\mathbb{Z}) \xrightarrow{\cup} H^2(G, \mathbb{Z}/q\mathbb{Z})\,.$$

From the exact sequence $0 \to \mathbb{Z}/q\mathbb{Z} \xrightarrow{q} \mathbb{Z}/q^2\mathbb{Z} \to \mathbb{Z}/q\mathbb{Z} \to 0$, we get the Bockstein homomorphism

$$B : H^1(G, \mathbb{Z}/q\mathbb{Z}) \longrightarrow H^2(G, \mathbb{Z}/q\mathbb{Z})$$

which is surjective and $\Delta$-invariant.

For a pro-$p$-group $P$ we denote by $P^i$, $i \geq 1$, the descending $q$-central series, i.e.

$$P^1 = P \quad \text{and} \quad P^{i+1} = (P^i)^q[P^i, P] \quad \text{for } i \geq 1\,.$$

Let

$$1 \longrightarrow F \longrightarrow \mathcal{F} \longrightarrow \Delta \longrightarrow 1$$

be an exact sequence of profinite groups where $F$ is a finitely generated pro-$p$-group. Obviously, $G^i$ and $F^i$ are normal open subgroups of $\mathcal{G}$ and $\mathcal{F}$ respectively.

**Proposition 1.3** *With the notation as above let $q > 2$ and $m \geq 2$. Assume that there exists a surjection*

$$\varphi_{m+1} : \mathcal{G} \qquad \mathcal{F}/F^{m+1}\,.$$

*Then there exists a surjection*

$$\varphi : \mathcal{G} \qquad \mathcal{F}$$

*inducing the surjection* $\varphi_m : \mathcal{G} \overset{\varphi_{m+1}}{\quad} \mathcal{F}/F^{m+1} \overset{can}{\quad} \mathcal{F}/F^m\,.$

**Proof:** Assume that we have already found a surjection

$$\varphi_{i+1} : \mathcal{G} \qquad \mathcal{F}/F^{i+1}$$

for $i \geq m$ which induces $\varphi_m$, and let $\varphi_i : \mathcal{G} \overset{\varphi_{i+1}}{\quad} \mathcal{F}/F^{i+1} \overset{can}{\quad} \mathcal{F}/F^i$.

Let $\gamma, x_0, \ldots, x_n$ be a minimal system of generators of $G$ such that $x_k \in \ker(\chi)$ for $k \geq 0$ and $\chi(\gamma) = 1 - q$.

*Claim*: The group $F^{i+1}/F^{i+2}$ is generated by elements of the form

$$w^q[w, \bar{\gamma}] \bmod F^{i+2}\,, \quad [w, \bar{x}_k] \bmod F^{i+2},\ k \geq 0, \quad w \in F^i\,,$$

where $\bar{\gamma}, \bar{x}_k \in F$ are lifts of $\varphi_2(\gamma), \varphi_2(x_k) \in F/F^2$.

This shown in [3] prop. 5(i) (observe, that we have a surjection $G/G^{i+1} \twoheadrightarrow F/F^{i+1}$, and so the group $F/F^{i+1}$ is generated by the elements $\bar{\gamma}, \bar{x}_k \bmod F^{i+1}$).

Consider the diagram with exact line

$$(*) \qquad \begin{array}{ccccccccc} & & & & & & \mathcal{G} & & \\ & & & & & & \varphi_i & & \\ 1 & & F^i/F^{i+2} & & \mathcal{F}/F^{i+2} & & \mathcal{F}/F^i & & 1\,. \end{array}$$

Since $i \geq m \geq 2$, the group $F^i/F^{i+2}$ is abelian,

$$[F^i, F^i] \subseteq F^{2i} \subseteq F^{i+2}\,,$$

and we consider $F^i/F^{i+2}$ as a $\mathcal{G}$-module via $\varphi_i$. Since the dualizing module $I$ of $G$ is isomorphic to $\mathbb{Q}_p/\mathbb{Z}_p$ as an abelian group, the canonical exact sequence

$$0 \longrightarrow F^{i+1}/F^{i+2} \longrightarrow F^i/F^{i+2} \longrightarrow F^i/F^{i+1} \longrightarrow 0$$

induces a $\Delta$-invariant exact sequence

$$0 \longrightarrow \mathrm{Hom}_G(F^i/F^{i+1}, I) \longrightarrow \mathrm{Hom}_G(F^i/F^{i+2}, I) \longrightarrow \mathrm{Hom}_G(F^{i+1}/F^{i+2}, I)\,.$$

Let $f \in \mathrm{Hom}_G(F^i/F^{i+2}, I)$. Then

$$\begin{array}{rcl} f([w, \bar{x}_k] \bmod F^{i+2}) & = & f(w \bmod F^{i+2})^{x_k - 1} = 0 \quad \text{for } k \geq 0, \\ f(w^q[w, \bar{\gamma}] \bmod F^{i+2}) & = & f(w \bmod F^{i+2})q + f(w \bmod F^{i+2})^{\gamma - 1} \\ & = & f(w \bmod F^{i+2})(q - q) = 0\,. \end{array}$$

Using the claim, we see that $f$ vanishes on $F^{i+1}/F^{i+2}$, and so

$$\mathrm{Hom}_G(F^i/F^{i+1}, I) \xrightarrow{\sim} \mathrm{Hom}_G(F^i/F^{i+2}, I)\,.$$

By duality, cf. [5] (3.7.6), (3.7.1), (3.4.6), we get

$$H^2(G, F^i/F^{i+2}) \xrightarrow{\sim} H^2(G, F^i/F^{i+1})\,,$$

and so

$$H^2(G, F^i/F^{i+2})^{\Delta} \xrightarrow{\sim} H^2(G, F^i/F^{i+1})^{\Delta}\,.$$

Since the order of $\Delta$ is prime to $p$, we obtain the isomorphism

$$H^2(\mathcal{G}, F^i/F^{i+2}) \xrightarrow{\sim} H^2(\mathcal{G}, F^i/F^{i+1})\,.$$

Now we prove that the embedding problem $(*)$ is solvable. For this we have to show that the 2-class

$$[\beta_i] \in H^2(\mathcal{F}/F^i, F^i/F^{i+2})$$

is mapped to zero under the inflation map $inf = \varphi_i^*$,

$$H^2(\mathcal{F}/F^i, F^i/F^{i+2}) \xrightarrow{inf} H^2(\mathcal{G}, F^i/F^{i+2})\,,$$

where $\beta_i$ is the 2-cocycle corresponding to the group extension in $(*)$, see [5] (9.4.2). From the commutative exact diagram

$$\begin{array}{ccccccccc}
 & & & & & & \mathcal{G} & & \\
 & & & & & \overset{\varphi_{i+1}}{\swarrow} & \downarrow \varphi_i & & \\
1 & \longrightarrow & F^i/F^{i+2} & \longrightarrow & \mathcal{F}/F^{i+2} & \longrightarrow & \mathcal{F}/F^i & \longrightarrow & 1 \quad \beta_i \\
 & & \downarrow can & & \downarrow & & \| & & \\
1 & \longrightarrow & F^i/F^{i+1} & \longrightarrow & \mathcal{F}/F^{i+1} & \longrightarrow & \mathcal{F}/F^i & \longrightarrow & 1 \quad \alpha_i
\end{array}$$

we get a commutative diagram

$$\begin{array}{ccc}
H^2(\mathcal{F}/F^i, F^i/F^{i+1}) & \xrightarrow{\varphi_i^*} & H^2(\mathcal{G}, F^i/F^{i+1}) \\
\uparrow can_* & & \uparrow can_* \\
H^2(\mathcal{F}/F^i, F^i/F^{i+2}) & \xrightarrow{\varphi_i^*} & H^2(\mathcal{G}, F^i/F^{i+2})\,.
\end{array}$$

Since there exists the solution $\varphi_{i+1}$ for the embedding problem $\alpha_i$, we have $\varphi_i^*([\alpha_i]) = 0$, and so

$$can_* \circ \varphi_i^*([\beta_i]) = \varphi_i^* \circ can_*([\beta_i]) = \varphi_i^*([\alpha_i]) = 0\,.$$

From the injectivity of the map $can_*$ on the right-hand side of the diagram above it follows that $\varphi_i^*([\beta_i]) = 0$, and so there exists a solution

$$\varphi_{i+2} : \mathcal{G} \longrightarrow \mathcal{F}/F^{i+2}$$

of the embedding problem corresponding to $\beta_i$. This homomorphism is necessarily surjective and induces $\varphi_m$, because $\varphi_i$ has these properties, cf. [5] (3.9.2).

Using a compactness argument, we get in the limit a surjection $\varphi : \mathcal{G} \twoheadrightarrow \mathcal{F}$ inducing $\varphi_m$. This finishes the proof of the proposition. □

In the following let $p$ be an odd prime number and let $\Delta = \langle\sigma\rangle \cong \mathbb{Z}/2\mathbb{Z}$ be cyclic of order 2. We denote, as usual, the $(\pm)$-eigenspaces of a $\mathbb{Z}_p[\Delta]$-module $M$ by $M^{\pm}$.

**Proposition 1.4** *Let $p$ be an odd prime number and let $G$ be a $p$-Demuškin group of rank $n+2$, $n \geq 0$, with dualizing module $I$ and invariant $q = \#(I^G) < \infty$. Assume that $\Delta \cong \mathbb{Z}/2\mathbb{Z}$ acts on $G$. Then the following holds:*

(i) *If $H^2(G, \mathbb{Z}/p\mathbb{Z}) = H^2(G, \mathbb{Z}/p\mathbb{Z})^-$, then $G_\Delta$ is a free pro-$p$-group of rank $n/2+1$.*

(ii) *If $H^2(G, \mathbb{Z}/p\mathbb{Z}) = H^2(G, \mathbb{Z}/p\mathbb{Z})^+$, then $G_\Delta$ is a $p$-Demuškin group of rank $m+2$, $0 \leq m \leq n$, with invariant $q$ and dualizing module $I$.*

**Proof:** We start with the following remark. Since $\mathrm{Aut}(I) \cong \mathbb{Z}_p^\times$ is abelian, the surjection $G \twoheadrightarrow \Gamma$ factors through $G_\Delta$. With the notation of the proof of proposition 1.1, it follows that $N \cap G$ has infinite index in $G$ and therefore $cd_p(N) = cd_p(N \cap G) \leq 1$, cf. [5] chap.II, §7, ex.3. Using the Hochschild-Serre spectral sequence and the fact that $\mathrm{Hom}(N, \mathbb{Z}/p\mathbb{Z}) = 0$, we see that $inf_2$ is an isomorphism, and so

$$H^2(G_\Delta) \cong H^2(G)^\Delta.$$

(i) By proposition 1.1, $G_\Delta$ is a free pro-$p$-group. Since the non-degenerate pairing

$$H^1(G) \times H^1(G) \xrightarrow{\cup} H^2(G) \cong \mathbb{Z}/p\mathbb{Z}$$

is $\Delta$-invariant, it follows that

$$\dim_{\mathbb{F}_p} H^1(G)^{\pm} = n/2 + 1\,.$$

Therefore

$$\dim_{\mathbb{F}_p} H^1(G_\Delta) = \dim_{\mathbb{F}_p} H^1(G)^\Delta = n/2 + 1\,.$$

(ii) If $H^2(G) = H^2(G)^+$, then $H^2(G_\Delta) \cong H^2(G)$, and we obtain a non-degenerate pairing

$$H^1(G_\Delta) \times H^1(G_\Delta) \xrightarrow{\cup} H^2(G_\Delta) \cong \mathbb{Z}/p\mathbb{Z}$$

showing that $G_\Delta$ is a $p$-Demuškin group. Finally, since $G_\Delta$ is non-trivial and its rank has to be even, it follows that $\dim_{\mathbb{F}_p} H^1(G_\Delta) \geq 2$, and since $\ker(G \twoheadrightarrow G_\Delta)$ acts trivially on $I$, we have $\#(I^{G_\Delta}) = \#(I^G) = q$ and $I$ is also the dualizing module of $G_\Delta$. □

## 2 Free Quotients

As before, let $G$ be a $p$-Demuškin group of rank $n+2$ with dualizing module $I$ and assume that $2 < q < \infty$. We are interested in quotients of $G$ which are free pro-$p$-groups. First we calculate the possible ranks of such quotients.

**Proposition 2.1** *Let $G$ be a Demuškin group of rank $n+2$ with finite invariant $q > 2$ and let $F$ be a free quotient of $G$. Then*

(i) *$H^1(F, \mathbb{Z}/q\mathbb{Z})$ lies in the kernel of the Bockstein homomorphism and*

(ii) *$H^1(F, \mathbb{Z}/q\mathbb{Z})$ is a totally isotropic free $\mathbb{Z}/q\mathbb{Z}$-submodule of $H^1(G, \mathbb{Z}/q\mathbb{Z})$ with respect to the pairing given by the cup-product.*

*In particular,*

$$\operatorname{rank} F \leq \frac{n}{2} + 1 \,.$$

**Proof:** Since $F$ is free, $H^1(F, \mathbb{Z}/q\mathbb{Z})$ is a free $\mathbb{Z}/q\mathbb{Z}$-module. The commutative diagram

$$\begin{array}{ccc} H^1(G, \mathbb{Z}/q\mathbb{Z}) \times H^1(G, \mathbb{Z}/q\mathbb{Z}) & \overset{\cup}{\longrightarrow} & H^2(G, \mathbb{Z}/q\mathbb{Z}) \\ \uparrow {\scriptstyle (inf\,, inf\,)} & & \uparrow {\scriptstyle inf} \\ H^1(F, \mathbb{Z}/q\mathbb{Z}) \times H^1(F, \mathbb{Z}/q\mathbb{Z}) & \overset{\cup}{\longrightarrow} & H^2(F, \mathbb{Z}/q\mathbb{Z}) = 0 \end{array}$$

shows that $H^1(F, \mathbb{Z}/q\mathbb{Z})$ is a totally isotropic $\mathbb{Z}/q\mathbb{Z}$-submodule of $H^1(G, \mathbb{Z}/q\mathbb{Z})$, and so $\dim_{\mathbb{F}_p} H^1(F, \mathbb{Z}/p\mathbb{Z}) = \operatorname{rank}_{\mathbb{Z}/q\mathbb{Z}} H^1(F, \mathbb{Z}/q\mathbb{Z}) \leq n/2 + 1$. From the commutative diagram

$$\begin{array}{ccc} H^1(G, \mathbb{Z}/q\mathbb{Z}) & \overset{B}{\longrightarrow} & H^2(G, \mathbb{Z}/q\mathbb{Z}) \\ \uparrow {\scriptstyle inf} & & \uparrow {\scriptstyle inf} \\ H^1(F, \mathbb{Z}/q\mathbb{Z}) & \overset{B}{\longrightarrow} & H^2(F, \mathbb{Z}/q\mathbb{Z}) = 0 \end{array}$$

follows that $H^1(F, \mathbb{Z}/q\mathbb{Z}) \subseteq \ker(B)$. □

Recall that $\Gamma$ is the canonical quotient $G/\ker(\chi)$ of $G$, where $\chi : G \longrightarrow \operatorname{Aut}(I)$ is the character given by the action of $G$ on $I$, i.e.

$$\chi_0 : \Gamma \longrightarrow \operatorname{Aut}(I) \,.$$

**Lemma 2.2** *The submodules $H^1(\Gamma, \mathbb{Z}/q\mathbb{Z})$ and $\ker B$ of $H^1(G, \mathbb{Z}/q\mathbb{Z})$ are orthogonal to each other, more precisely*

$$H^1(\Gamma, \mathbb{Z}/q\mathbb{Z}) = (\ker B)^{\perp} \,.$$

**Proof:** Consider the commutative diagram of non-degenerate pairings

$$\begin{array}{ccccccc} H^1(G, {}_qI) & \times & H^1(G, \mathbb{Z}/q\mathbb{Z}) & \cup & H^2(G, {}_qI) & & \mathbb{Z}/q\mathbb{Z} \\ \delta & & B & & & & \\ H^0(G, {}_qI) & \times & H^2(G, \mathbb{Z}/q\mathbb{Z}) & \cup & H^2(G, {}_qI) & & \mathbb{Z}/q\mathbb{Z} \end{array}$$

which is induced by the exact sequences

$$0 \longrightarrow \mathbb{Z}/q\mathbb{Z} \xrightarrow{q} \mathbb{Z}/q^2\mathbb{Z} \longrightarrow \mathbb{Z}/q\mathbb{Z} \longrightarrow 0 \quad \text{and} \quad 0 \longrightarrow {}_qI \xrightarrow{q} {}_{q^2}I \longrightarrow {}_qI \longrightarrow 0\,.$$

Let $\zeta_{q^2}$ be a generator of ${}_{q^2}I$, and so $\zeta_q = (\zeta_{q^2})^q$ is a generator of ${}_qI$. By definition of the homomorphism $\delta$, i.e.

$$\delta : H^0(G, {}_qI) \longrightarrow H^1(G, {}_qI), \quad (\zeta_q)^i \;\mapsto \Big\{ G \to {}_qI,\; g \mapsto ((\zeta_{q^2})^i)^{g-1} \Big\}\,,$$

we see that the image of $\delta$ is $H^1(\Gamma, {}_qI)$. Since the pairings above are non-degenerated, it follows that $H^1(\Gamma, \mathbb{Z}/q\mathbb{Z}) = H^1(\Gamma, {}_qI)$ is orthogonal to $\ker B$, and since $\mathrm{rank}_{\mathbb{Z}/q\mathbb{Z}}(\ker B)^\perp = 1$, we prove the lemma. $\square$

**Proposition 2.3** *Let $G$ be a $p$-Demuškin group of rank $n+2$ with finite invariant $q > 2$ and let $F$ be a free factor of $G$ of rank $n/2+1$. Then the canonical surjection $G \twoheadrightarrow \Gamma$ factors through $F$, i.e. there is a commutative diagram*

$$\begin{array}{ccc} G & & \Gamma \\ & F\,. & \end{array}$$

**Proof:** Suppose the contrary. Then there exists an open subgroup $G'$ of $G$ which has a surjection

$$(G')^{ab} \twoheadrightarrow (F')^{ab} \times \Gamma'\,,$$

where $F'$ is the image of $G'$ in $F$ under the projection $G \twoheadrightarrow F$ and $\Gamma'$ is the image of $G'$ under the projection $G \twoheadrightarrow \Gamma$. Let $q' = \#(I^{G'}) = \#(I^{\Gamma'})$. Since $F'$ is free, it follows that $H^1(F', \mathbb{Z}/q'\mathbb{Z})$ is a totally isotropic submodule of $H^1(G', \mathbb{Z}/q'\mathbb{Z})$ and contained in $\ker B'$ by proposition 2.1. From lemma 2.2 we know that $H^1(\Gamma', \mathbb{Z}/q'\mathbb{Z})$ is orthogonal to $\ker B'$, and so also to $H^1(F', \mathbb{Z}/q'\mathbb{Z})$. Thus $H^1(F', \mathbb{Z}/q'\mathbb{Z}) \oplus H^1(\Gamma', \mathbb{Z}/q'\mathbb{Z})$ is totally isotropic. But $H^1(F', \mathbb{Z}/q'\mathbb{Z})$ is a maximal totally isotropic $\mathbb{Z}/q'\mathbb{Z}$-submodule of $H^1(G', \mathbb{Z}/q'\mathbb{Z})$ of rank $d \cdot n/2 + 1$, where $d = (G : G')$. This contradiction proves the proposition. $\square$

For the existence of free quotients of Demuškin groups we have the following

**Theorem 2.4** *Let $G$ be a $p$-Demuškin group of rank $n+2$ with finite invariant $q > 2$ and let $\varDelta$ be a finite abelian group of exponent $p-1$ acting on $G$. Let $V$ be a $\mathbb{Z}/q\mathbb{Z}$-submodule of $H^1(G, \mathbb{Z}/q\mathbb{Z})$ such that*

(i) *$V$ is $\mathbb{Z}/q\mathbb{Z}$-free and $\varDelta$-invariant,*

(ii) *$V$ is totally isotropic with respect to the pairing given by the cup-product,*

(iii) *$V$ lies in the kernel of the Bockstein map $B : H^1(G, \mathbb{Z}/q\mathbb{Z}) \to H^2(G, \mathbb{Z}/q\mathbb{Z})$.*

*Then there exists a $\varDelta$-invariant surjection*

$$G \qquad F$$

*onto a free quotient $F$ of $G$ such that $H^1(F, \mathbb{Z}/q\mathbb{Z}) = V$.*

**Proof:** Let

$$1 \longrightarrow R \longrightarrow F_{n+2} \longrightarrow G \longrightarrow 1$$

be a minimal presentation of $G$, where $F_{n+2}$ is a free pro-$p$-group of rank $n+2$. Using lemma 1.2, we extend the action of $\varDelta$ to $F_{n+2}$. Let $\gamma, x_0, \ldots, x_n$ be a basis of $F_{n+2}$ such that

(i) each element of the basis of $F_{n+2}$ generates a $\varDelta$-invariant subgroup isomorphic to $\mathbb{Z}_p$ on which $\varDelta$ acts by some character $\psi : \varDelta \to \mu_{p-1}$,

(ii) $R$, as a normal subgroup of $F_{n+2}$, is generated by the element

$$w = (x_0)^q [x_0, \gamma][x_1, x_2][x_3, x_4] \cdot \cdots \cdot [x_{n-1}, x_n] \cdot f$$

where $f \in (F_{n+2})^3$,

(iii) $V^* = \mathrm{Hom}(V, \mathbb{Z}/q\mathbb{Z})$ has a basis $\{v_i \mod (F_{n+2})^2, 1 \le i \le r = \mathrm{rank}_{\mathbb{Z}/q\mathbb{Z}} V\}$ such that

$$\{v_1, \ldots, v_r\} \quad \text{is a subset of} \quad \{\gamma, x_1, \ldots, x_n\}$$

and, if $v_i = x_{j(i)}$, then $x_{j(i)+1} \notin \{v_1, \ldots, v_r\}$ or $x_{j(i)-1} \notin \{v_1, \ldots, v_r\}$ according to whether $j(i)$ is odd or even.

Such a basis exists: by [1] prop.(2.3), we find a basis of $F_{n+2}$ with the property (i) (since $\varDelta$ is abelian of exponent $p-1$, the $\mathbb{Z}/q\mathbb{Z}[\varDelta]$-module $F_{n+2}/(F_{n+2})^2$ decomposes into a sum of $\mathbb{Z}/q\mathbb{Z}[\varDelta]$-modules which are free $\mathbb{Z}/q\mathbb{Z}$-modules of rank equal to 1). The $\varDelta$-invariance of the cup-product and the Bockstein homomorphism implies that we find a basis satisfying (i) and (ii), cf. [8] lemma 3 in the case where $\varDelta \cong \mathbb{Z}/2\mathbb{Z}$. Using the assumptions on $V$, we can also satisfy (iii).

Let $N$ be the normal subgroup of $F_{n+2}$ generated by the set

$$\{\gamma, x_k, \, 0 \le k \le n\} \setminus \{v_1, \ldots, v_r\},$$

then $F := F_{n+2}/N$ is a free pro-$p$-group of rank $r$, $N$ is $\Delta$-invariant and we have $R \subseteq N(F_{n+2})^3$ by the properties (ii) and (iii) of the basis $\gamma, x_0, \ldots, x_n$. Thus the $\Delta$-invariant surjection

$$F_{n+2} \twoheadrightarrow F/F^3 = F_{n+2}/N(F_{n+2})^3$$

factors through a $\Delta$-invariant surjection

$$G \qquad F/F^3 \, .$$

Applying proposition 1.3, we get a $\Delta$-invariant surjection from $G$ onto a free pro-$p$-group $F$ which induces a surjection $G \twoheadrightarrow F/F^2 \cong F_{n+2}/N(F_{n+2})^2$.

By construction, we have $F/F^2 \cong V^*$, and so $H^1(F, \mathbb{Z}/q\mathbb{Z}) = V$. This finishes the proof of the theorem. □

Now we consider free quotients of a Demuškin group $G$ which are invariant under a given $\Delta$-action of $G$, where $\Delta$ is a group of order 2.

**Corollary 2.5** *Let $p$ be an odd prime number and let $G$ be a $p$-Demuškin group of rank $n+2$, $n \geq 0$, with finite invariant $q$. Let $\Delta \cong \mathbb{Z}/2\mathbb{Z}$ acting on $G$ such that $H^2(G, \mathbb{Z}/q\mathbb{Z}) = H^2(G, \mathbb{Z}/q\mathbb{Z})^-$. Let*

$$u^+,\, u^- \geq 0 \quad \textit{be integers such that} \quad u^+ + u^- = n/2 \, .$$

*Then there exists a $\Delta$-invariant surjection*

$$\varphi : G \qquad F$$

*such that*

(i) *$F$ is a free pro-$p$-group of rank $n/2 + 1$,*

(ii) $\mathrm{rank}_{\mathbb{Z}_p}(F^{ab})^+ = u^+ + 1$ *and* $\mathrm{rank}_{\mathbb{Z}_p}(F^{ab})^- = u^-$.

**Proof:** Since $H^2(G, \mathbb{Z}/q\mathbb{Z}) = H^2(G, \mathbb{Z}/q\mathbb{Z})^-$, the submodules $H^1(G, \mathbb{Z}/q\mathbb{Z})^\pm$ are maximal totally isotropic with respect to the cup-product pairing, and so

$$\mathrm{rank}_{\mathbb{Z}/q\mathbb{Z}} H^1(G, \mathbb{Z}/q\mathbb{Z})^\pm = n/2 + 1.$$

Let

$$V = V^+ \oplus V^-$$

where $V^+$ is a free $\mathbb{Z}/q\mathbb{Z}$-submodule of $H^1(G, \mathbb{Z}/q\mathbb{Z})^+$ of rank $1 + u^+$ containing $H^1(\Gamma, \mathbb{Z}/q\mathbb{Z})$, and $V^-$ is a free $\mathbb{Z}/q\mathbb{Z}$-submodule of $\ker B^-$ of rank $u^-$ being orthogonal to $V^+$, i.e. $V$ is maximal totally isotropic and contained in $\ker B$.

It is easy to see that $V^-$ exists (see also the remark below): by lemma 2.2 $H^1(\Gamma,\mathbb{Z}/q\mathbb{Z}) \subseteq (\ker B^-)^\perp$, and since

$$\mathrm{rank}_{\mathbb{Z}/q\mathbb{Z}} H^1(G,\mathbb{Z}/q\mathbb{Z})^+ - \mathrm{rank}_{\mathbb{Z}/q\mathbb{Z}} V^+ = n/2+1-(1+u^+) = u^- ,$$

there exists a free submodule of $\ker B^-$ of rank $u^-$ which is orthogonal to $V^+$. Now theorem 2.4 gives us a free $\Delta$-invariant quotient $F$ of $G$ of rank $n/2+1$ such that

$$H^1(F,\mathbb{Z}/q\mathbb{Z}) = V \cong (\mathbb{Z}/q\mathbb{Z}[\Delta]^+)^{u^++1} \oplus (\mathbb{Z}/q\mathbb{Z}[\Delta]^-)^{u^-}.$$

Since $F^{ab}$ is free $\mathbb{Z}_p$-module, we obtain assertion (ii). □

**Remark:** Explicitly, we get a submodule $V$ with the properties as above in the following way: let

$$1 \longrightarrow R \longrightarrow F_{n+2} \longrightarrow G \longrightarrow 1$$

be a minimal presentation of $G$, where $F_{n+2}$ is a free pro-$p$-group of rank $n+2$ with the extended action of $\Delta$. Let $\gamma, x_0, \dots, x_n$ be a basis of $F_{n+2}$ such that $R$ is generated by the element

$$w = (x_0)^q [x_0,\gamma][x_1,x_2][x_3,x_4]\cdot\dots\cdot[x_{n-1},x_n]\cdot f ,$$

$f \in (F_{n+2})^3$, and

$$\begin{array}{lll} \gamma^\sigma = \gamma\cdot a , & x_i^\sigma = x_i \cdot a_i & \text{for } i=2,4,\dots,n , \\ x_0^\sigma = x_0^{-1}\cdot b, & x_i^\sigma = x_i^{-1}\cdot b_i & \text{for } i=1,3,5,\dots,n-1 . \end{array}$$

with $a$, $b$, $a_i$, $b_i \in (F_{n+2})^2$. Such a basis exists by the $\Delta$-invariance of the cup-product and the Bockstein homomorphism, cf. [8] lemma 3. If we put

$$\begin{array}{lll} \gamma' := \gamma\cdot a^{\frac{1}{2}} , & x_i' := x_i \cdot a_i^{\frac{1}{2}} & \text{for } i=2,4,\dots,n , \\ x_0' := b^{-\frac{1}{2}}\cdot x_0 , & x_i' := b_i^{-\frac{1}{2}}\cdot x_i & \text{for } i=1,3,5,\dots,n-1 , \end{array}$$

then

$$\begin{array}{lll} (\gamma')^\sigma = \gamma' , & (x_i')^\sigma = x_i' & \text{for } i\geq 2 \text{ even,} \\ (x_0')^\sigma = (x_0')^{-1}, & (x_i')^\sigma = (x_i')^{-1} & \text{for } i\geq 1 \text{ odd,} \end{array}$$

and

$$w = (x_0')^q [x_0',\gamma'][x_1',x_2'][x_3',x_4']\cdot\dots\cdot[x_{n-1}',x_n']\cdot f'$$

where $f' \in (F_{n+2})^3$. Let $u = 2u^+ - 1$. If we denote $x \bmod F^2$ by $\bar{x}$, then the dual of

$$\begin{aligned} V^* &:= \mathbb{Z}/q\mathbb{Z}\cdot\bar{\gamma} \oplus \bigoplus_{i=1,3,\dots,u} \mathbb{Z}/q\mathbb{Z}\cdot\bar{x}_{i+1} \oplus \bigoplus_{i=u+3,\dots,n} \mathbb{Z}/q\mathbb{Z}\cdot\bar{x}_{i-1} \\ &\cong (\mathbb{Z}/q\mathbb{Z}[\Delta]^+)^{u^++1} \oplus (\mathbb{Z}/q\mathbb{Z}[\Delta]^-)^{u^-} \end{aligned}$$

gives an example for a submodule with the properties (i)-(iii) in the proof of corollary 2.5. The free quotient of $G$ is obtained in the following way: if

$$N' = (x'_0, \underbrace{x'_1, x'_3, \ldots, x'_u}_{u^+\text{-times}}, \underbrace{x'_{u+3}, \ldots, x'_n}_{u^-\text{-times}}) \trianglelefteq F_{n+2}\,,$$

then $F = F_{n+2}/N'$ is a free pro-$p$-group of rank $n/2+1$, $N'$ is $\Delta$-invariant, $R \subseteq N(F_{n+2})^3$ and $V^* = F/F^2$. Using proposition 1.3 we get the desired quotient of $G$.

With the notation and assumptions of corollary 2.5, we make for a $\Delta$-invariant free quotient $F$ of $G$ of rank $n/2+1$ the following

**Definition 2.6** *We call the tuple* $(u^+, u^-)$ *the* **signature** *of* $F$*, if*

$$F/F^2 \cong (\mathbb{Z}/q\mathbb{Z}[\Delta]^+)^{u^+ +1} \oplus (\mathbb{Z}/q\mathbb{Z}[\Delta]^-)^{u^-}\,.$$

It is easy to see that the signature of a maximal free quotient $F$ of $G$ does not determine $F$, but if the signature is equal to $(n/2, 0)$, then we have the following proposition.

**Proposition 2.7** *Let $p$ be an odd prime number and let $\Delta$ be of order 2. Let $G$ be a $p$-Demuškin group of rank $n+2$ with finite invariant $q$ on which $\Delta$ acts such that $H^2(G, \mathbb{Z}/p\mathbb{Z})^\Delta = 0$. Let $F$ be a free $\Delta$-invariant quotient of $G$ of rank $n/2+1$, i.e. the canonical surjection*

$$G \qquad F$$

*is $\Delta$-invariant. If the induced action of $\Delta$ on $F/F^2$ is trivial, i.e. $F$ has signature $(n/2, 0)$, then $F$ is equal to the maximal quotient $G_\Delta$ of $G$ with trivial $\Delta$-action. In particular, a free quotient of $G$ with the properties above is unique.*

**Proof:** As in the remark after the proof of corollary 2.5, we find generators of $F$ on which $\Delta$ acts trivially, and so $F$ has a trivial $\Delta$-action. Thus we have a surjection

$$\varphi : G_\Delta \qquad F\,.$$

Since $G_\Delta$ is free of rank $n/2+1 = \dim_{\mathbb{F}_p} H^1(F, \mathbb{Z}/p\mathbb{Z})$ by proposition 1.4(i), it follows that $\varphi$ is an isomorphism. Thus $F$ is the maximal quotient of $G$ with trivial $\Delta$-action. □

Mathematisches Institut
der Universität Heidelberg
Im Neuenheimer Feld 288
69120 Heidelberg
Germany

e-mail: wingberg@mathi.uni-heidelberg.de